\documentclass[11pt]{amsart}

\usepackage{amsmath}
\usepackage{amssymb}
\usepackage{mathrsfs}
\usepackage[all]{xy}

\newtheorem{theorem}{Theorem}[section]
\newtheorem{claim}[theorem]{Claim}

\theoremstyle{definition}
\newtheorem{definition}[theorem]{Definition}

\newtheorem{question}[theorem]{Question}

\theoremstyle{remark}

\newcount\skewfactor
\def\mathunderaccent#1#2 {\let\theaccent#1\skewfactor#2
\mathpalette\putaccentunder}
\def\putaccentunder#1#2{\oalign{$#1#2$\crcr\hidewidth
\vbox to.2ex{\hbox{$#1\skew\skewfactor\theaccent{}$}\vss}\hidewidth}}
\def\name{\mathunderaccent\tilde-3 }

\def\smallbox#1{\leavevmode\thinspace\hbox{\vrule\vtop{\vbox
   {\hrule\kern1pt\hbox{\vphantom{\tt/}\thinspace{\tt#1}\thinspace}}
   \kern1pt\hrule}\vrule}\thinspace}

\newcommand{\cf}{{\rm cf}}

\def\qedref#1{$\qed_{\reforiginal{#1}}$}

\title{Aristotelian poetry}
\author{Shimon Garti}
\address{Einstein Institute of Mathematics,
 The Hebrew University of Jerusalem,
 Jerusalem 91904, Israel}
\email{shimon.garty@mail.huji.ac.il}

\author{Saharon Shelah}
\address{Institute of Mathematics
 The Hebrew University of Jerusalem,
 Jerusalem 91904, Israel
 and  Department of Mathematics
 Rutgers University
 New Brunswick, NJ 08854, USA}
\email{shelah@math.huji.ac.il}
\urladdr{http://www.math.rutgers.edu/\char`\~shelah}
\thanks{The research was supported by Israel Science Foundation Grant no. 1838/19. This is publication 1219a of the second author}


\begin{document}
\let\labeloriginal\label
\let\reforiginal\ref
\def\ref#1{\reforiginal{#1}}
\def\label#1{\labeloriginal{#1}}

\begin{abstract}
Jing Zhang proved in \cite{MR4094551} the consistency of $\binom{\omega_2}{\omega_1}\rightarrow \binom{n}{\omega_1}_\omega$ for every $n\in\omega$ with the negative relation $\binom{\omega_2}{\omega_1}\nrightarrow \binom{\omega}{\omega_1}_\omega$.
We improve the consistency strength and obtain this relation from an $\omega_1$-Erd\H{o}s cardinal.
\end{abstract}

\maketitle

\newpage

\section{Introduction}

We consider a problem about polarized relations with infinitely many colors.
Recall that $\binom{\alpha}{\beta}\rightarrow\binom{\gamma}{\delta}_\chi$ denotes the statement that for every coloring $c:\alpha\times\beta\rightarrow\chi$ one can find $A\subseteq\alpha,B\subseteq\beta$ and $i\in\chi$ such that ${\rm otp}(A)=\gamma,{\rm otp}(B)=\delta$ and $c''(A\times B)=\{i\}$.
The most investigated case is $\beta=\kappa$ and $\alpha=\kappa^+$ for some infinite cardinal $\kappa$.
Our convention is that $\beta\leq\alpha$ and we tacitly assume that $\chi\in\kappa$ to avoid trivialities.

In this paper, $\chi$ will be an infinite cardinal, a fact which implies almost immediately some consistency strength with respect to positive relations.
One can see this from the following (unpublished) result of Galvin.

\begin{claim}
\label{conclmgalvin}
If there exists a Kurepa tree then $\binom{\omega_2}{\omega_1}\nrightarrow\binom{2}{\omega_1}_\omega$.
\end{claim}

\par\noindent\emph{Proof}. \newline
Suppose that $f,g:\omega_1\rightarrow\omega$.
We shall say that $f$ and $g$ are \emph{almost disjoint} iff $|\{\beta\in\omega_1:f(\beta)=g(\beta)\}|\leq\aleph_0$.
Let $\mathscr{T}$ be a Kurepa tree.
We intend to build a collection $\mathscr{F} = \{f_\alpha:\alpha\in\omega_2\}\subseteq{}^{\omega_1}\omega$ such that $\mathscr{F}$ is a family of pairwise almost disjoint functions.

For every $\beta\in\omega_1$ fix an enumeration $(t_{\beta n}:n\in\omega)$ of the elements of $\mathcal{L}_\beta(\mathscr{T})$.
Let $(b_\alpha:\alpha\in\omega_2)$ be an enumeration of the $\omega_1$-branches of $\mathscr{T}$ (or some of them, if there are more than $\aleph_2$ branches).
For every $\alpha\in\omega_2$ define $f_\alpha:\omega_1\rightarrow\omega$ as follows:
$$
f_\alpha(\beta)=m \Leftrightarrow b_\alpha\cap\mathcal{L}_\beta(\mathscr{T})=t_{\beta m}.
$$
Notice that if $\alpha_0<\alpha_1<\omega_2$ then for some $\beta_0\in\omega_1$ we have $b_{\alpha_0}\upharpoonright\beta_0 = b_{\alpha_1}\upharpoonright\beta_0$ and $b_{\alpha_0}(\beta)\neq b_{\alpha_1}(\beta)$ whenever $\beta\in[\beta_0,\omega_1)$.
By the definition of our functions we see that $f_{\alpha_0}(\beta)\neq f_{\alpha_1}(\beta)$ for every $\beta\in[\beta_0,\omega_1)$.

Letting $\mathscr{F}=\{f_\alpha:\alpha\in\omega_2\}$ we may conclude that $\mathscr{F}$ is almost disjoint.
All we need now is to convert such a family to a coloring which exemplifies the negative relation $\binom{\omega_2}{\omega_1}\nrightarrow\binom{2}{\omega_1}_\omega$.
To this end, define a coloring $c:\omega_2\times\omega_1\rightarrow\omega$ as follows:
$$
c(\alpha,\beta)=f_\alpha(\beta).
$$
If $\alpha_0<\alpha_1<\omega_2$ and $B\in[\omega_1]^{\omega_1}$ then for a sufficiently large $\beta\in B$ we have $c(\alpha_0,\beta)=f_{\alpha_0}(\beta)\neq f_{\alpha_1}(\beta)=c(\alpha_1,\beta)$, so we are done.

\hfill \qedref{conclmgalvin}

Silver proved that the non-existence of Kurepa trees is equiconsistent with the existence of a strongly inaccessible cardinal.
From Silver's theorem we deduce that the positive relation $\binom{\omega_2}{\omega_1}\rightarrow \binom{2}{\omega_1}_\omega$ requires at least a strongly inaccessible cardinal.
But the consistency strength of this statement is stronger, and it is located around an $\omega_1$-Erd\H{o}s cardinal as proved by Donder and Levinski in \cite{MR1024901}.

In the present work we are interested in a stepping-up phenomenon, reflected in the following theorem of Baumgartner who showed in \cite{MR1034564} that if $\binom{\omega_2}{\omega_1}\rightarrow \binom{2}{\omega_1}_\omega$ then $\binom{\omega_2}{\omega_1}\rightarrow \binom{n}{\omega_1}_\omega$ for every $n\in\omega$.
As recorded in \cite{MR1034564}, this statement was also proved by Donder and Levinski.
A natural question is how far can this stepping-up theorem reach.

\begin{question}
\label{q111} Does the positive relation $\binom{\omega_2}{\omega_1}\rightarrow \binom{2}{\omega_1}_\omega$ imply $\binom{\omega_2}{\omega_1}\rightarrow \binom{\omega}{\omega_1}_\omega$?
\end{question}

This is \cite[Question 1.11]{MR4101445}, and our goal is to give a negative answer.
Namely, we shall force the positive relation $\binom{\omega_2}{\omega_1}\rightarrow \binom{2}{\omega_1}_\omega$ (and hence a positive relation will hold with sets of size $n$ for every natural number $n$), but simultaneously we will force $\binom{\omega_2}{\omega_1}\nrightarrow \binom{\omega}{\omega_1}_\omega$.

This result can be viewed as a reflection of an old debate between Plato and Aristotle, see \cite{aristo}.
According to the common interpretation, Aristotle denied the concept of \emph{actual infinity} but accepted the idea of \emph{potential infinity}.
Therefore, one can think of an arbitrarily large natural number $n$ but the collection of all natural numbers does not exist.
The main result of this paper is a combinatorial version of this attitude.
During the refereeing process we learned that the main result of the paper has been proved by Jing Zhang in \cite{MR4094551}.
We indicate, however, that our proof reduces the consistency strength of the statement to an $\omega_1$-Erd\H{o}s cardinal.

Our notation is mostly standard and follows, e.g., the conventions of \cite{MR506523}.
However, we employ the Jerusalem forcing notation, so $p\leq{q}$ means that the condition $p$ is weaker than the condition $q$.
For background in the partition calculus we suggest \cite{MR2768681}, \cite{MR3075383} and \cite{MR795592}.

\newpage

\section{Polarized relations with infinitely many colors}

An infinite cardinal $\lambda$ is called $\omega_1$-Erd\H{o}s iff $\lambda\rightarrow(\omega_1)^{<\omega}_2$.
These cardinals were studied combinatorially in \cite{MR0202613}.
Several deep results were established by Silver in his dissertation with respect to these cardinals.
For more on this subject we direct the reader to Silver's paper \cite{MR409188}.

\begin{theorem}
\label{thmmt} Assume that $\lambda$ is $\omega_1$-Erd\H{o}s.
Then one can force $\binom{\omega_2}{\omega_1}\rightarrow \binom{n}{\omega_1}_\omega$ for every $n\in\omega$ and simultaneously $\binom{\omega_2}{\omega_1}\nrightarrow \binom{\omega}{\omega_1}_\omega$.
\end{theorem}

\par\noindent\emph{Proof}. \newline
We may assume that Martin's axiom with $2^\omega=\omega_2$ holds in the ground model, since one can force it while preserving the properties of $\lambda$.
We define a forcing notion $\mathbb{Q}$.
The conditions in $\mathbb{Q}$ will approximate a coloring $c:\lambda\times\omega_1\rightarrow\omega$.
This coloring will be a witness (in the generic extension) to the negative relation $\binom{\lambda}{\omega_1}\nrightarrow\binom{\omega}{\omega_1}_\omega$.
Concomitantly, $\mathbb{Q}$ will collapse the cardinals below $\lambda$ to $\aleph_1$ but $\lambda$ will be preserved.
Thus $\lambda=\omega_2$ in the generic extension and the desired negative relation will be established.
The positive relation $\binom{\omega_2}{\omega_1}\rightarrow \binom{n}{\omega_1}_\omega$ will be a reminiscent of the fact that $\lambda$ is $\omega_1$-Erd\H{o}s in the ground model.

Let us define the forcing conditions and the order.
A condition $p\in\mathbb{Q}$ is a quintuple $(\varepsilon,U,f,g,\mathcal{A})=(\varepsilon^p,U^p,f^p,g^p,\mathcal{A}^p)$ such that:
\begin{enumerate}
\item [$(a)$] $\varepsilon\in\omega_1$.
\item [$(b)$] $U\subseteq\lambda-\omega_1,|U|\leq\aleph_1$.
\item [$(c)$] $f$ and $g$ are functions whose domain is $U$.
\item [$(d)$] If $\alpha\in{U}$ then $f(\alpha)$ is a one-to-one function from $\varepsilon$ into $U\cup\alpha$.
\item [$(e)$] If $\alpha\in{U}$ then $g(\alpha)$ is a function from $\varepsilon$ into $\omega$.
\item [$(f)$] $\mathcal{A}\subseteq[U]^{\aleph_0}$ and $|\mathcal{A}|\leq\aleph_1$.
\end{enumerate}
We shall refer to $\varepsilon^p$ as \emph{the height} of the condition $p$.
Suppose that $p,q\in\mathbb{Q}$.
We shall say that $p\leq_{\mathbb{Q}}q$ iff $\varepsilon^p\leq\varepsilon^q, U^p\subseteq U^q,\mathcal{A}^p\subseteq\mathcal{A}^q$ and the following requirements are met:
\begin{enumerate}
\item [$(\aleph)$] If $\alpha\in{U^p}$ then $f^p(\alpha)\subseteq f^q(\alpha)$ and $g^p(\alpha)\subseteq g^q(\alpha)$.
\item [$(\beth)$] If $u\in{\mathcal{A}^p}$ and $\varepsilon^p\leq\varepsilon<\varepsilon^q$ then there are $\alpha,\beta\in{u}$ such that $\alpha\neq\beta$ and $g^q(\alpha)(\varepsilon)\neq g^q(\beta)(\varepsilon)$.
\end{enumerate}
Recall that $f(\alpha)$ and $g(\alpha)$ are functions, so the inclusion in the order definition is inclusion of functions.

The \emph{working part} of the condition is $g$, which approximates a witness to the desired negative relation.
The r\^{o}le of $f$ is to collapse all the cardinals strictly below $\lambda$ to $\aleph_1$.
One can verify that $\leq_{\mathbb{Q}}$ is a partial ordering, thus $\mathbb{Q}$ is a forcing notion.

If $(p_n:n\in\omega)$ is an increasing sequence of conditions in $\mathbb{Q}$ then it has an upper bound (actually, a least upper bound).
To see this, let $\varepsilon=\bigcup_{n\in\omega}\varepsilon^{p_n}$ so $\varepsilon\in\omega_1$.
Let $U=\bigcup_{n\in\omega}U^{p_n}$ so $U\subseteq\lambda-\omega_1$ and $|U|\leq\aleph_1$.
Similarly, let $\mathcal{A}=\bigcup_{n\in\omega}\mathcal{A}^{p_n}$ and for every $\alpha\in{U}$ let $f(\alpha)=\bigcup\{f^{p_n}(\alpha):\alpha\in U^{p_n},n\in\omega\}$ and $g(\alpha)=\bigcup\{g^{p_n}(\alpha):\alpha\in U^{p_n},n\in\omega\}$.
Set $q=(\varepsilon,U,f,g,\mathcal{A})$ and notice that $p_n\leq_{\mathbb{Q}}q$ for every $n\in\omega$.

Fix $G\subseteq\mathbb{Q}$ which is $V$-generic.
It follows that $\aleph_1$ is preserved in $V[G]$.
Since all the objects in a given condition are of size at most $\aleph_1$ and $\lambda^{\omega_1}=\lambda$ one can deduce from a Delta-system argument that $\mathbb{Q}$ is $\lambda$-cc, thus $\lambda$ is preserved in $V[G]$ as well.
Notice, however, that if $\kappa\in(\omega_1,\lambda)$ then the $f$-part of the conditions will add a mapping from $\omega_1$ onto $\kappa$, thus $\kappa$ will be an ordinal of size $\aleph_1$ in $V[G]$.
One concludes that $\lambda=\omega_2^{V[G]}$ and this fact will be proved anon in a formal way.

Let us define two types of dense open subsets of $\mathbb{Q}$.
Firstly, for every $\alpha<\beta<\lambda$ and each $\varepsilon\in\omega_1$ we define:
\[D_{\alpha\beta\varepsilon}=\{p\in\mathbb{Q}:\alpha,\beta\in U^p, \varepsilon\leq\varepsilon^p,\beta\in{\rm rang}(f^q(\alpha))\}.\]
Let us show that each $D_{\alpha\beta\varepsilon}$ is open and dense.
The openness of $D_{\alpha\beta\varepsilon}$ follows directly from the definitions.
For density, suppose that $p\in\mathbb{Q}$ but $p\notin D_{\alpha\beta\varepsilon}$.
Let us try to define a condition $q$ so that $p\leq{q}\in D_{\alpha\beta\varepsilon}$.

Choose $\varepsilon^q\geq\varepsilon$ such that $\varepsilon^p+\omega<\varepsilon^q$.
Define $U^q=U^p\cup\{\alpha,\beta\}\cup\omega_1$.
For every $\gamma\in U^p$ let $f^q(\gamma)$ be a one-to-one function from $\varepsilon^q$ into $U^q$ such that $f^q(\gamma)\supseteq f^p(\gamma)$.
This is possible since $\varepsilon^q\in\omega_1$ and $|U^q|=\aleph_1$.
Likewise, for $\gamma=\beta$ we make sure that $\alpha\in{\rm rang}(f^q(\beta))$.
This is possible since $\varepsilon^p<\varepsilon^q$ and hence one can pick $i\in\varepsilon^q-\varepsilon^p$ and set $f^q(\beta)(i)=\alpha$.

The interesting mission is $g^q$, and here we have to be attentive to part $(\beth)$ of the order definition.
We define $\mathcal{A}^q=\mathcal{A}^p$ and if $\gamma\in U^p$ then $g^q(\gamma)\upharpoonright\varepsilon^p = g^p(\gamma)\upharpoonright\varepsilon^p$, and if $\gamma\in U^q-U^p$ then $g^q(\gamma)\upharpoonright\varepsilon^p$ is constantly zero.
Now for every $\gamma\in U^q$ and every $\zeta\in[\varepsilon^p,\varepsilon^q)$ we choose $g^q(\gamma)(\zeta)$ such that for every $u\in\mathcal{A}^p=\mathcal{A}^q$ the sequence $(g^q(\gamma)(\zeta):\gamma\in{u})$ is not constant.
This is possible due to Martin's axiom and $2^\omega=\omega_2$ (forced at the beginning of the proof), but we make the comment that adding $\omega_2$-many Cohen reals instead of MA will have the same effect.
Now $q\in D_{\alpha\beta\varepsilon}$ (this is the easy task) and $p\leq_{\mathbb{Q}}q$ (this is the challenging part) so $D_{\alpha\beta\varepsilon}$ is dense.

For every $\alpha\in[\omega_1,\lambda)$ let $\name{g}_\alpha=\bigcup\{g^p(\alpha):p\in G,\alpha\in U^p\}$ and let $\name{f}_\alpha=\bigcup\{f^p(\alpha):p\in G,\alpha\in U^p\}$.
From the density of the $D_{\alpha\beta\varepsilon}$s we deduce that $\name{f}_\alpha$ maps $\omega_1$ onto $\alpha$ and hence every $\kappa\in(\omega_1,\lambda)$ is collapsed to $\aleph_1$.
Likewise, $\name{g}_\alpha$ is a function from $\omega_1$ to $\omega$.

Secondly, for every $u\in[\lambda]^{\aleph_0}$ define:
\[E_u=\{p\in\mathbb{Q}:u\in\mathcal{A}^p\}.\]
Clearly, $E_u$ is open.
For the density of $E_u$ notice that if $p\notin E_u$ then one can choose $U^q\supseteq U^p$ such that $u\subseteq U^q$ and then one can add $u$ to $\mathcal{A}^q$ while defining $f^q$ and $g^q$ according to the requirements of $\leq_{\mathbb{Q}}$.

To accomplish the proof of the negative relation we combine the $\name{g}_\alpha$s into a (name of) a single coloring $\name{c}:(\lambda-\omega_1)\times\omega_1\rightarrow\omega$.
Thus for $\alpha\in\lambda-\omega_1,\beta\in\omega_1$ let $\name{c}(\alpha,\beta)=\name{g}_\alpha(\beta)$ and let $c=\name{c}[G]$.
By the above considerations, $c$ is a coloring from $(\omega_2-\omega_1)\times\omega_1$ into $\omega$ in $V[G]$.

Suppose that $u\in[\omega_2-\omega_1]^{\aleph_0}$ and $v\in[\omega_1]^{\aleph_1}$.
Let $p\in G$ be any condition which forces this fact.
Choose $\alpha,\beta\in{u}$ and $\varepsilon\in v$ such that $\alpha\neq\beta$ and $\varepsilon>\varepsilon^p$.
Let $q$ be so that $p\leq q\in D_{\alpha\beta\varepsilon}$.
Now $q\Vdash\name{c}(\alpha,\varepsilon) = \name{g}_\alpha(\varepsilon) = \name{g}(\alpha)(\varepsilon) \neq \name{g}(\beta)(\varepsilon) = \name{g}_\beta(\varepsilon) = \name{c}(\beta,\varepsilon)$ and hence the negative relation $\binom{\omega_2}{\omega_1}\nrightarrow \binom{\omega}{\omega_1}_\omega$ is established.

It remains to prove that the positive relation $\binom{\omega_2}{\omega_1}\rightarrow \binom{n}{\omega_1}_\omega$ holds in $V[G]$ whenever $n\in\omega$.
By \cite{MR1034564} it is sufficient to concentrate on the case of $n=2$.
Recall that $\lambda\rightarrow(\omega_1)^{<\omega}_2$, and Martin's axiom with $2^\omega=\omega_2$ holds in the ground model.

Let $\chi=\lambda^+$ and let $<^*$ be a well-ordering of $\mathcal{H}(\chi)$.
Suppose that $p\in\mathbb{Q}$ and $p\Vdash\name{c}:\omega_2\times\omega_1\rightarrow\omega$.
By the assumption that $\lambda$ is $\omega_1$-Erd\H{o}s we choose a set $U\subseteq\lambda$ such that ${\rm otp}(U)=\omega_1$ and $U$ is indiscernible in the structure $(\mathcal{H}(\chi),\in,<^*,p,\mathbb{Q},c)$.
That is, if $\alpha,\beta\in{U}$ and $\alpha<\beta$ then $U-\beta$ is indiscernible over $\alpha=\{\gamma:\gamma\in\alpha\}$.

For our argument we would like to recast the whole story to a relatively small elementary submodel.
Let $N$ be the Skolem hull of $U\cup\{p,\mathbb{Q},c\}$ in $(\mathcal{H}(\chi),\in,<^*)$.
Notice that ${\rm otp}(N\cap\lambda)=\omega_1, |N\cap\omega_1|=\aleph_0$ and $p,\mathbb{Q},c\in{N}$.
Denote $N\cap\omega_1$ by $\delta$, so $\delta$ is a countable ordinal.
Let $\mathbb{R}=\mathbb{Q}\upharpoonright{N}$.
Applying at $\mathbb{R}$ the argument for the $\lambda$-cc of $\mathbb{Q}$ one concludes that $\mathbb{R}$ is $\aleph_1$-cc.

Let $G\subseteq\mathbb{R}$ be a directed set which intersects every dense subset $D$ of $\mathbb{Q}$ which satisfies $D\in{N}$.
The existence of $G$ is based on the assumption that Martin's axiom holds in the ground model, and the fact the the collection of dense sets has the finite intersection property.
Let $q\in\mathbb{Q}$ be an upper bound of $G$.
Choose a $V$-generic $H\subseteq\mathbb{Q}$ so that $q\in{H}$.

We define a function $d:N\cap\lambda\rightarrow\omega$ by $d(\alpha)=c(\alpha,\delta)$.
Since $\Vdash_{\mathbb{Q}}\omega_1^V=\omega_1$ and since $U\subseteq{N}\cap\lambda,|U|=\aleph_1$ there is a natural number $m\in\omega$ in $V[H]$ such that $W=\{\alpha\in{U}:d(\alpha)=m\}\in[U]^{\aleph_1}$.
For $\{\alpha_0,\alpha_1\}\subseteq{W}$ let $A_{\alpha_0\alpha_1}=\{\zeta\in\omega_1: c(\alpha_0,\zeta)=c(\alpha_1,\zeta)=m\}$.
Notice that $|A_{\alpha_0\alpha_1}|=\aleph_1$ whenever $\{\alpha_0,\alpha_1\}\subseteq{W}$.
To see this, suppose that $\{\alpha_0,\alpha_1\}\subseteq{W}$ and $\name{\zeta}$ is a $\mathbb{Q}$-name for $\sup(\name{A}_{\alpha_0\alpha_1})$.
Notice that $\name{\zeta}\in{N}$ (it is definable in $N$) and hence $\name{\zeta}[H]\in{N}$.
By the definition of $W$ and $d$ we see that $\delta\in\name{A}_{\alpha_0\alpha_1}$ so necessarily $\name{\zeta}[H]=\omega_1$ and hence $A_{\alpha_0\alpha_1}=\name{A}_{\alpha_0\alpha_1}[H]$ is of size $\aleph_1$.
This means that $c\upharpoonright(\{\alpha_0,\alpha_1\}\times A_{\alpha_0\alpha_1})$ is constantly $m$, so the proof is accomplished.

\hfill \qedref{thmmt}

One concludes from the above theorem that the stepping-up argument of Baumgartner breaks down at $\aleph_0$.
But the breaking point corresponds to the domain of the coloring.
Thus if $\kappa\geq 2^\omega$ then $\binom{\kappa^{++}}{\kappa^+}\rightarrow\binom{2}{\kappa^+}_\omega$ implies $\binom{\kappa^{++}}{\kappa^+}\rightarrow\binom{\aleph_0}{\kappa^+}_\omega$.
Nevertheless, for successor cardinals the breaking point is strictly smaller than the small component of the polarized relation.
If $\kappa$ is strongly inaccessible then the stepping-up argument goes up to $\kappa$.

\begin{theorem}
\label{thmstepupinac} If $\kappa$ is strongly inaccessible and $\binom{\kappa^+}{\kappa}\rightarrow\binom{2}{\kappa}_{<\kappa}$ then $\binom{\kappa^+}{\kappa}\rightarrow\binom{\gamma}{\kappa}_{<\kappa}$ for every $\gamma\in\kappa$.
\end{theorem}

\par\noindent\emph{Proof}. \newline
Suppose that $\gamma\in\kappa$ and $\binom{\kappa^+}{\kappa}\nrightarrow\binom{\gamma}{\kappa}_{<\kappa}$.
Let $\chi$ be such that $\binom{\kappa^+}{\kappa}\nrightarrow\binom{\gamma}{\kappa}_{\chi}$.
By monotonicity we may assume that $\chi=|\gamma|^+$.
Our goal is to prove that $\binom{\kappa^+}{\kappa}\nrightarrow\binom{2}{\kappa}_{\chi}$, thus proving our statement.

Let $f:\kappa^+\times\kappa\rightarrow\chi$ be a witness to the negative relation $\binom{\kappa^+}{\kappa}\nrightarrow\binom{\gamma}{\kappa}_{\chi}$.
For every $\xi\in\kappa^+$ let $f_\xi=f\upharpoonright\{\xi\}\times\kappa$, so $f_\xi\in{}^\kappa\chi$.
Assume that $x\in[\kappa^+]^\gamma$.
Let $\tau=\tau(x)\in\kappa$ be the first ordinal for which if $\beta\in[\tau,\kappa)$ then $|\{f_\xi(\beta):\xi\in x\}|\geq 2$.
Such an ordinal exists since $f$ witnesses $\binom{\kappa^+}{\kappa}\nrightarrow\binom{\gamma}{\kappa}_{\chi}$.

For every $\xi\in\kappa^+$ let $h_\xi:\kappa\rightarrow\xi$ be onto.
By induction on $\xi\in\kappa^+$ we would like to define $g_\xi:\kappa\rightarrow\chi\times\chi$, so fix $\xi\in\kappa^+$ and assume that $g_\eta$ is already defined for every $\eta<\xi$.
We will assume, further, that each $g_\eta(\zeta)$ is a pair of the form $(f_\eta(\zeta),j)$ for some $j\in\chi$.
Set:
$$
C_\xi=\{\alpha\in\kappa:\forall x\in[h''_\xi \alpha]^\gamma,\tau(x\cup\{\xi\})<\alpha\}.
$$
Clearly, $C_\xi$ is a club subset of $\kappa$.
For every $\alpha\in\kappa$ let $\delta(\alpha)=\sup(C_\xi\cap\alpha)$, so $\delta(\alpha)\leq\alpha$ and $\delta(\alpha)\in C_\xi$.
For every $\alpha\in\kappa$ define:
$$y_\alpha=\{\eta\in h''_\xi \delta(\alpha): f_\eta(\alpha)=f_\xi(\alpha)\}.$$
Notice that $|y_\alpha|\leq|\gamma|<\chi$.

Choose $j\in\chi$ so that $\eta\in y_\alpha\Rightarrow j\neq j_\eta$, stipulating $g_\eta(\alpha)=(f_\eta(\alpha),j_\eta)$, and let $g_\xi(\alpha)=(f_\xi(\alpha),j)$.
This defines a function $g_\xi:\kappa\rightarrow\chi\times\chi$.
Define $g:\kappa^+\times\kappa\rightarrow\chi\times\chi$ by $g(\xi,\alpha)=g_\xi(\alpha)$.

Suppose that $\eta<\xi<\kappa^+$.
Choose $\delta\in C_\xi-(\eta+1)$.
For every $\alpha\geq\delta$ one has $g_\xi(\alpha)=g(\xi,\alpha)\neq g(\eta,\alpha)=g_\eta(\alpha)$, either because $f_\xi(\alpha)\neq f_\eta(\alpha)$ or because $j^\alpha_\xi\neq j^\alpha_\eta$.
This implies that $\binom{\kappa^+}{\kappa}\nrightarrow\binom{2}{\kappa}_{\chi}$ as required.

\hfill \qedref{thmstepupinac}

The following seems to be open:

\begin{question}
\label{qstronglyinac} Is it consistent that $\kappa$ is a strongly inaccessible cardinal and yet $\binom{\kappa^+}{\kappa}\nrightarrow \binom{2}{\kappa}_{<\kappa}$?
\end{question}

If $\kappa$ is weakly compact then it is easy to see that $\binom{\kappa^+}{\kappa}\rightarrow \binom{2}{\kappa}_{<\kappa}$, and hence $\binom{\kappa^+}{\kappa}\rightarrow\binom{\gamma}{\kappa}_{<\kappa}$ for every $\gamma\in\kappa$.
However, if one requires a stationary set in the small component then a negative result can be forced over weakly compact cardinals.

Hajnal proved in \cite{MR281605} that if $\kappa$ is measurable then $\binom{\kappa^+}{\kappa}\rightarrow \binom{\kappa}{\kappa}_{<\kappa}$ and a slight modification gives the relation $\binom{\kappa^+}{\kappa}\rightarrow \binom{\tau}{\kappa}_{<\kappa}$ for every $\tau\in\kappa^+$.
The result reappeared in a paper of Chudnovskii who claimed that it holds at every weakly compact cardinal as well.
He did not supply the argument, and it is still opaque whether the statement is true or not.
The best known results are $\binom{\kappa^+}{\kappa}\rightarrow \binom{\kappa^n}{\kappa}_{<\kappa}$ for every $n\in\omega$ and $\binom{\kappa^+}{\kappa}\rightarrow \binom{\tau}{\kappa}_m$ for every $\tau\in\kappa^+,m\in\omega$, due to Jones in \cite{MR2275863}.
The latter is an improvement upon a result of Wolfsdorf from \cite{MR603336}.

Thus the first open problem is whether $\binom{\kappa^+}{\kappa}\rightarrow \binom{\kappa^\omega}{\kappa}_\omega$ holds at every weakly compact cardinal, see \cite[Question 8.3]{MR2768681}.
Another way to phrase this question is by asking whether one can distinguish measurability and weak compactness on the ground of the polarized relation.
We suggest below a possible direction by considering a stationary set in the small component of the monochromatic product.

Let us say that $\binom{\lambda}{\kappa}\rightarrow\binom{\tau}{stat}_\chi$ if for every $c:\lambda\times\kappa\rightarrow\chi$ one can find $A\subseteq\lambda,{\rm otp}(A)=\tau$ and a stationary set $B\subseteq\kappa$ such that $c\upharpoonright(A\times B)$ is constant.
Similar expressions should be interpreted accordingly, e.g. if $\mathscr{U}$ is an ultrafilter over $\kappa$ then $\binom{\lambda}{\kappa}\rightarrow\binom{\tau}{\mathscr{U}}_\chi$ means that the small component $B\subseteq\kappa$ satisfies $B\in\mathscr{U}$.

As we shall see, if $\kappa$ is measurable, $\tau\in\kappa^+$ and $\mathscr{U}$ is a normal ultrafilter over $\kappa$ then $\binom{\kappa^+}{\kappa}\rightarrow\binom{\tau}{\mathscr{U}}_{<\kappa}$ (this is a result of Kanamori, \cite{MR498138}) and since $\mathscr{U}$ is normal this means that $\binom{\kappa^+}{\kappa}\rightarrow\binom{\tau}{stat}_{<\kappa}$.
But if $\kappa$ is weakly compact then one can force the negation of similar statements.

In order to do it we shall use Kurepa trees, adapted to inaccessible cardinals.
If one imitates the classical definition of Kurepa trees then the tree ${}^{<\kappa}2$ is Kurepa whenever $\kappa$ is strongly inaccessible, thus the concept of Kurepa trees becomes uninteresting.
We shall use the traditional substitute:

\begin{definition}
\label{defslim} Slim Kurepa trees. \newline
A $\kappa$-tree $\mathscr{T}$ is a slim Kurepa tree iff $|\mathcal{L}_\beta(\mathscr{T})|\leq|\beta|$ for every $\beta\in\kappa$ and the cardinality of the set of cofinal branches of $\mathscr{T}$ is at least $\kappa^+$.
\end{definition}

Before proceeding let us exclude from the discussion a seemingly stronger concept.
Call $\mathscr{T}$ \emph{very slim} if $|\mathcal{L}_\beta(\mathscr{T})|<|\beta|$ for every $\beta\in\kappa$.
Let us show that the number of cofinal branches in such trees must be small.

\begin{claim}
\label{clmveryslim} Suppose that $\kappa=\cf(\kappa)>\aleph_0$.
Then there are no very slim $\kappa$-Kurepa trees.
\end{claim}

\par\noindent\emph{Proof}. \newline
Let $\mathscr{T}$ be a very slim $\kappa$-tree.
We intend to prove that the cardinality of the set of cofinal branches of $\mathscr{T}$ is strictly less than $\kappa$.
For every $\beta\in\kappa$ let $\theta_\beta=|\mathcal{L}_\beta(\mathscr{T})|$, so $\theta_\beta<|\beta|$.
By Fodor's lemma there are a stationary $S\subseteq\kappa$ and a fixed cardinal $\theta\in\kappa$ such that $\beta\in S\Rightarrow\theta_\beta=\theta$.
Choose such a pair $(S,\theta)$ where $\theta$ is minimal.
Let $\mathcal{B}$ be the set of cofinal branches of $\mathscr{T}$.
We claim that $|\mathcal{B}|\leq\theta$.

Assume towards contradiction that $|\mathcal{B}|>\theta$ and let $\{b_\alpha:\alpha\in\theta^+\}\subseteq\mathcal{B}$ be a set of pairwise distinct cofinal branches.
For every pair of ordinals $\{\alpha,\delta\}\in[\theta^+]^2$ let $\beta_{\alpha\delta}\in\kappa$ be so that $b_\alpha(\beta_{\alpha\delta})\neq b_\delta(\beta_{\alpha\delta})$.
Let $\xi=\bigcup\{\beta_{\alpha\delta}:\{\alpha,\delta\}\in[\theta]^2\}$.
Notice that $\xi\in\kappa$ since $\kappa=\cf(\kappa)>\theta^+$.
Choose $\beta\in S$ such that $\beta>\xi$.
For every $\alpha\in\theta^+$ let $x_\alpha=b_\alpha\cap\mathcal{L}_\beta(\mathscr{T})$.
By the above choices if $\alpha\neq\delta$ then $x_\alpha\neq x_\delta$.
However, $|\mathcal{L}_\beta(\mathscr{T})|=\theta<\theta^+$, so this is impossible.

\hfill \qedref{clmveryslim}

It is easy to see that if $\kappa$ is measurable (in fact, ineffability is sufficient) then there are no slim $\kappa$-Kurepa trees.
On the other hand, if $\kappa$ is weakly compact then one can force such a tree over $\kappa$ while preserving weak compactness.
In the following theorem we shall see that the existence of such a tree gives a negative polarized relation in the sense of stationary sets.

\begin{theorem}
\label{thmwcneg} Let $\kappa$ be weakly compact and let $S=S^\kappa_\omega\cap{\rm Card}$.
Then one can force the existence of a coloring $c:\kappa^+\times S\rightarrow\omega$ with no monochromatic product $A\times T$ such that ${\rm otp}(A)=\kappa$ and $T$ is a stationary subset of $S$.
\end{theorem}

\par\noindent\emph{Proof}. \newline
Let $\mathscr{T}$ be a slim Kurepa tree over $\kappa$, and let $\{b_\alpha:\alpha\in\kappa^+\}$ be a set of (pairwise distinct) cofinal branches in $\mathscr{T}$.
By the above claim we may assume that $|\mathcal{L}_\beta(\mathscr{T})|=|\beta|$ for every $\beta\in\kappa$.

For each $\beta\in S$ we choose a decomposition of $\mathcal{L}_\beta(\mathscr{T})$ of the form $(I^\beta_n:n\in\omega)$ such that $|I^\beta_n|=\theta^\beta_n<\beta$.
We define $c:\kappa^+\times{S}\rightarrow\omega$ by $c(\alpha,\beta)=n$ iff $b_\alpha\cap\mathcal{L}_\beta(\mathscr{T})\in I^\beta_n$.
We claim that our statement is witnessed by the coloring $c$.
Indeed, suppose that $A\subseteq\kappa^+$ with ${\rm otp}(A)=\kappa$ and $T$ is a stationary subset of $S$, for which $c\upharpoonright(A\times{T})$ is constantly $n$.
For every $\beta\in T$ one has $\theta^\beta_n<\beta$, so by Fodor's lemma there is a stationary $T_0\subseteq{T}$ and a fixed $\theta\in\kappa$ such that $\beta\in T_0\Rightarrow \theta^\beta_n=\theta$.

Let $A_0$ consist of the first $\theta^+$ elements of $A$.
For every $\alpha,\alpha'\in A_0$ there is $\gamma(\alpha,\alpha')\in\kappa$ such that $b_\alpha(\gamma)\neq b_{\alpha'}(\gamma)$ whenever $\gamma\geq\gamma(\alpha,\alpha')$.
Since $|A_0|=\theta^+$ one can fix $\gamma_0\in\kappa$ such that $\gamma(\alpha,\alpha')<\gamma_0$ for every pair of ordinals $\alpha$ and $\alpha'$ from $A_0$.
Choose $\beta\in T_0$ such that $\beta\geq\gamma_0$.
For every $\alpha\in A_0$ let $x_\alpha=b_\alpha\cap\mathcal{L}_\beta(\mathscr{T})\in I^\beta_n$.
Notice that if $\alpha\neq\alpha'$ then $x_\alpha\neq x_{\alpha'}$ since $\beta\geq\gamma_0>\gamma(\alpha,\alpha')$.
Let $W=\{x_\alpha:\alpha\in A_0\}$.
It follows that $|W|=\theta^+$.
However, $W\subseteq I^\beta_n$ and $|I^\beta_n|=\theta$, a contradiction.

\hfill \qedref{thmwcneg}

The positive result of Jones at weakly compact cardinals, namely $\binom{\kappa^+}{\kappa}\rightarrow \binom{\kappa^n}{\kappa}_{<\kappa}$ for every $n\in\omega$, is limited essentially by the Milner-Rado paradox.
Hence a reasonable attempt to obtain $\binom{\kappa^+}{\kappa}\nrightarrow \binom{\kappa^\omega}{\kappa}_\omega$ would be something based on the decomposition of $\kappa^\omega$ (or a larger ordinal) according to this paradox.
However, this cannot follow from the Milner-Rado paradox in a straightforward way, since the paradox applies to successors of measurable cardinals as well.

Therefore, one has to employ some combinatorial principle which (consistently) holds at weakly compact cardinals and fails at measurable cardinals.
Slim Kurepa trees are viable candidates as reflected in Theorem \ref{thmwcneg}.
Since ineffable cardinals do not carry such trees, one may ask:

\begin{question}
\label{qineffable} Does the positive relation $\binom{\kappa^+}{\kappa}\rightarrow \binom{\tau}{\kappa}_{<\kappa}$ holds for every $\tau\in\kappa^+$ when $\kappa$ is ineffable?
\end{question}

\newpage

\bibliographystyle{amsplain}
\bibliography{arlist}

\end{document}